\documentclass[12pt]{article}
\usepackage{amsmath}
\usepackage{mathrsfs}
\usepackage{bbm}
\usepackage{amssymb,a4}
\usepackage{amsfonts,amssymb,amsthm}
\usepackage{extarrows}
\usepackage{color}
\usepackage[all,pdf]{xy}

\allowdisplaybreaks  

\newcommand{\pf}[1]{\noindent{\bf Proof.}#1\hfill{}$\Box$}

\newcommand{\Z}{{\mathbb Z}}
\newcommand{\R}{{\mathbb R}}
\newcommand{\C}{{\mathbb C}}
\newcommand{\al}{\alpha}
\newcommand{\be}{\beta}
\newcommand{\dt}{\delta}
\newcommand{\lmd}{\lambda}
\newcommand{\gm}{\gamma}
\newcommand{\vf}[1]{\varphi(#1)}
\newcommand{\e}{{\epsilon}}
\newcommand{\sgm}{\sigma}

\newcommand{\V}{\mathfrak V}
\newcommand{\W}{\mathfrak W}
\newcommand{\ing}{\in G}
\newcommand{\f}[2]{f\left(#1,#2\right)}
\newcommand{\g}[2]{g\left(#1,#2\right)}
\newcommand{\phif}[2]{\phi\left(#1,#2\right)}
\newcommand{\etaf}[1]{\eta(#1)}

\newtheorem{thm}{Theorem}[section]
\newtheorem{prop}[thm]{Proposition}
\newtheorem{lem}[thm]{Lemma}
\newtheorem{cor}[thm]{Corollary}

\theoremstyle{remark}

\begin{document}

\begin{center}
{\Large {\bf Compatible left-symmetric algebraic structures
              on high rank Witt and Virasoro algebras}}\\
\vspace{0.5cm}
\end{center}

\begin{center}
{Chengkang Xu\footnote{
The author is supported by the National Natural Science Foundation of China (No.11626157, 11801375),
the Science and Technology Foundation of Education Department of
Jiangxi Province (No. GJJ161044).}\\
Shangrao Normal University, Shangrao, Jiangxi, China\\
Email: xiaoxiongxu@126.com}
\end{center}

\begin{abstract}
We classify all graded compatible left-symmetric algebraic structures on high rank Witt algebras,
and classify all non-graded ones satisfying a minor condition.
Furthermore, graded compatible left-symmetric algebraic structures
on high rank Virasoro algebras are also classified.
\end{abstract}

\section{Introduction}

\def\theequation{1.\arabic{equation}}
\setcounter{equation}{0}

The left-symmetric algebra, also called the pre-Lie algebra by others,
has a long history back to Cayley in 1890's \cite{C}.
It is closely related to many branches of mathematics,
such as Lie algebra and Lie group, integrable system, Yang-Baxter equation and so on,
see \cite{Bor,DM1,DM2,LM} and the references therein.
An important way to study left-symmetric algebras is through their sub-adjacent Lie algebras.
In the finite dimension case, it is well known that
the sub-adjacent Lie algebra of a finite dimensional left-symmetric algebra
over a field of characteristic 0 is not semi-simple.
Therefore the classical representation theory
of finite dimensional semi-simple Lie algebras is useless here.
However, in the infinite dimension case,
a simple left-symmetric algebra can have a simple sub-adjacent Lie algebra.
The first example is the left-symmetric Witt algebra $Vec(n)$ \cite{Bur}.
Let $U$ be an associative commutative algebra and $\partial_1,\cdots,\partial_n$
commuting derivations of $U$.
Then the left-symmetric Witt algebra
$Vec(n)=\{\sum_{i=1}^nu_i\partial_i\mid u_i\in U\}$
has multiplication
$$(u\partial_i)\cdot(v\partial_j)=(u\partial_i(v))\partial_j,\ \ \ u,v\in U.$$
In particular, let $n=1, U=\C[t^{\pm1}]$ and $\partial=\frac{\partial}{\partial t}$,
one gets a simple left-symmetric algebra $W_1=\C[t^{\pm1}]\partial$,
whose sub-adjacent Lie algebra is the Witt algebra.
The Witt algebra $W_1$ is a simple complex Lie algebra with a basis $\{L_a\mid a\in\Z\}$
satisfying
\begin{equation}\label{eq1.1}
 [L_a,L_b]=(b-a)L_{a+b},\ \ \forall a,b\in\Z,
\end{equation}
whose universal central extension is the famous Virasoro algebra $Vir$,
which has commutation relation
\begin{equation}\label{eq1.2}
 [L_a,L_b]=(b-a)L_{a+b}+\frac K{12}(b^3-b)\dt_{a+b,0},\ \ \forall a,b\in\Z,
\end{equation}
where $K$ is a central element.

Motivated by this example, the authors in \cite{KCB} study
compatible left-symmetric algebraic structures (abbreviated as CLSAS) on the Witt algebra $W_1$
and the Virasoro algebra $Vir$.
They classified all graded CLSAS's on $W_1$,
namely the ones satisfying $L_aL_b=\f abL_{a+b}$ for some function $f$ on $\Z\times\Z$.
They proved that such a function $f$ must be
\begin{equation}\label{eq1.3}
 \f ab=\frac{(\al+b+a\al\theta)(1+\theta b)}{1+\theta(a+b)},\ \ \ \
 \al,\theta\in\C\text{ such that }\theta=0\text{ or }\theta^{-1}\not\in\Z,
\end{equation}
or
\begin{equation}\label{eq1.4}
 \f ab=\begin{cases}
  \lmd+b &\text{ if } \lmd+a+b\neq0;\\
  \frac{(\lmd+b)(\gm-\lmd-b)}{\gm-\lmd} &\text{ if } \lmd+a+b=0,
  \end{cases}\ \ \ \
  \gm\in\C,\lmd\in\Z\text{ such that }\gm\neq\lmd.
\end{equation}
Furthermore, graded CLSAS's on $Vir$ are also classified in \cite{KCB}.
They are given by
\begin{equation}\label{eq1.5}
 L_aL_b=\frac{b(1+\theta b)}{1+\theta(a+b)}L_{a+b}
  +\frac K{24}\left(b^3-b-(\theta-\theta^{-1})b^2\right)\dt_{a+b,0},
\end{equation}
where $K$ is a cental element and $\theta\in\C$ such that
$\theta^{-1}\not\in\Z, \mathrm{Re}\theta>0$
or $\mathrm{Re}\theta=0, \mathrm{Im}\theta>0$.
In \cite{TB}, and further in \cite{LGB},
non-graded CLSAS's on $W_1$ satisfying that
$L_aL_b=\f abL_{a+b}+\g abL_{a+b+\zeta}$ for some $\zeta\in\Z\setminus\{0\}$,
and functions $f,g$ on $\Z\times\Z$,
are classified.
They turn out to be Novikov algebras \cite{GD} and defined by
\begin{equation}\label{eq1.6}
 \f ab=\al+b,\  \g ab=\mu,\ \ \ \text{ where }\al\in\C, \mu\in\C\setminus\{0\}.
\end{equation}
Based on these classification results, the CLSAS's on
many other Lie algebras containing $W_1$ or $Vir$ as a subalgebra were considered.
For example, the $W$-algebra $W(2,2)$ \cite{CL1},
the Heisenberg-Virasoro algebra \cite{CL2}, and so on.

A common point of these Lie algebras is that they are all graded by the rank 1 free group $\Z$.
In this paper, we study CLSAS's on high rank Witt algebra $\W$
and high rank Virasoro algebra $\V$ \cite{PZ},
which are graded by a free subgroup $G$ of $\C$ of rank $\nu>1$.
The algebra $\W$ and $\V$ have same commutation relations as (\ref{eq1.1}) and (\ref{eq1.2})
with the index group $\Z$ replaced by $G$.
Our motivation also comes from the example $Vec(n)$.
Take $n=1, U=\C[t_1^{\pm1},\cdots,t_\nu^{\pm1}]$
the ring of Laurrent polynomials in $\nu$ variables, and
$$\partial_1=\sum_{i=1}^\nu\e_i t_i\frac{\partial}{\partial t_i},$$
where $\e_1,\cdots,\e_\nu$ is a $\Z$-basis of $G$,
then we get a CLSAS $(Vec(1),\cdot)$ on $\W$.
Moreover, if we take a nonzero $w\in U$ and define
$$(u\partial_1)\circ(v\partial_1)=(u\partial_1(v))\partial_1+(wuv)\partial_1,$$
then we get a non-graded CLSAS $(Vec(1),\circ)$ on $\W$.
Inspired by these two examples, it is natural to consider CLSAS's on $\W$ satisfying
\begin{equation}\label{eq1.7}
L_aL_b=\f abL_{a+b}+\g abL_{a+b+\zeta},
\end{equation}
where $\zeta\ing\setminus\{0\}$ and $f,g$ are functions on $G\times G$,
and furthermore to consider central extensions of these CLSAS's on $\W$
whose sub-adjacent Lie algebra is $\V$.
These CLSAS's on $\W$ and $\V$ turn out having the same expression as in (\ref{eq1.3})-(\ref{eq1.6})
with index group $\Z$ replaced by $G$.
As we all know, the algebra $Vir$ ($W_1$) and $\V$ (respectively $\W$)
bear many resemblances in both structure and representation theories.
Our study on the CLSAS's on $\W$ and $\V$ gives one more resemblance to these Lie algebras.

The paper is organized as follows.
In Section \ref{sec2}, we recall some notions and results about left-symmetric algebras,
$\W$ and $\V$.
Section \ref{sec3} is devoted to the classification of graded CLSAS's on $\W$,
and Section \ref{sec4} to the classification of non-graded CLSAS's on $\W$ satisfying (\ref{eq1.7}).
In the last section,
we compute central extensions of the graded left-symmetric algebras obtained in Section \ref{sec3},
whose sub-adjacent Lie algebra is $\V$.
Throughout this paper, the symbols $\Z,\Z_+,\R, \C$ refer to the set of integers, positive integers,
real numbers and complex numbers respectively.

\section{Left-symmetric algebras and high rank Virasoro algebras}\label{sec2}
\def\theequation{2.\arabic{equation}}
\setcounter{equation}{0}

In this section we recall the definition and properties of the left-symmetric algebra
and some basic results about high rank Witt and Virasoro algebras.

A {\em left-symmetric algebra} is a vector space $A$ over $\C$ equipped
with a bilinear product $(x,y)\mapsto xy$ satisfying that
\begin{equation}\label{eq2.1}
 (xy)z-x(yz)=(yx)z-y(xz)\text{ for any }x,y,z\in A.
\end{equation}
Set $[x,y]=xy-yx$.
This defines a Lie algebra structure on the vector space $A$,
called the {\em sub-adjacent Lie algebra} of $A$
and denoted by $L(A)$.
Meanwhile, $A$ is also called a {\em compatible left-symmetric algebraic structure}
on the Lie algebra $L(A)$.

We remark that there is a CLSAS on a Lie algebra $\mathcal{G}$
if and only if there exist a representation
$\rho:\mathcal G\longrightarrow\mathfrak{gl}(V)$ of $\mathcal G$,
and a bijective linear map $q:\mathcal G\longrightarrow V$ satisfying that
$$q[x,y]=\rho(x)q(y)-\rho(y)q(x),\ \ \ \forall x,y\in\mathcal G.$$

%

For later use we recall indecomposable modules $V$ over $\V$ with weight space decomposition
$V=\bigoplus_{a\ing}V_a$, where all $\dim V_a=1$ and
$V_a=\{v\in V\mid L_0v=(\lmd+a)v\}$ for some fixed $\lmd\in\C$.
There are three classes of such modules, $V_{\al,\be},A_\gm,B_\gm,$
with parameters $\al,\be,\gm\in\C$,
which share a same basis $\{v_a\mid a\ing\}$ and have $\V$-actions as follows.\\
$$\begin{aligned}
 V_{\al,\be}: &\  Kv_b=0,\ L_av_b=(\al+b+a\be)v_{a+b};\\
 A_\gm:       &\  Kv_b=0,\
                 L_av_b=\begin{cases}
                  (a+b)v_{a+b} &\text{ if } b\neq0,\\
                   a(\gm+a)v_a &\text{ if } b=0;
                 \end{cases}\\
 B_\gm:       &\  Kv_b=0,\
                 L_av_b=\begin{cases}
                   bv_{a+b} &\text{ if } a+b\neq0,\\
                   -a(\gm+a)v_0 &\text{ if } a+b=0.
                 \end{cases}
\end{aligned}$$
Clearly, these $V_{\al,\be},A_\gm,B_\gm$ are also indecomposable $\W$-modules.
\begin{thm}[\cite{Su}]\label{thm2.1}
 Any $\V$(or $\W$)-modules of intermediate series with all weight spaces being 1-dimensional is
 isomorphic to one of $V_{\al,\be},A_\gm,B_\gm$.
\end{thm}

\section{Graded CLSAS's on the high rank Witt algebra $\W$}
\label{sec3}
\def\theequation{3.\arabic{equation}}
\setcounter{equation}{0}

In this section we study the CLSAS on $\W$ with multipication
\begin{equation}\label{eq3.01}
 L_aL_b=f(a,b)L_{a+b}, \ \ \  \forall a,b\ing.
\end{equation}
We call such a CLSAS on $\W$ {\em $G$-graded}.
To avoid ambiguity, we denote such a CLSAS by $W$.
Then $W$ is a $\W$-module defined by the above equation.
Moreover, it follows from the definition of the left-symmetric algebra
and the Lie bracket of $\W$ that
\begin{align}
 &\f ab-\f ba=b-a,\label{eq3.02}\\
 &(b-a)\f{a+b}c=\f bc\f a{b+c}-\f ac\f b{a+c}.\label{eq3.03}
\end{align}

The next lemma is crucial for later computations.
\begin{lem}\label{lem3.1}
 $\f a0=\f00$ for any $a\ing$. Therefore $\f0a=\f00+a$.
\end{lem}
\pf{
Let $b=c=0$ in (\ref{eq3.03}) we get
$$-a\f a0=\f00\f a0-\f a0\f0a.$$
Apply (\ref{eq3.02}) and we obtain
\begin{equation}\label{eq3.04}
\f a0(\f a0-\f00)=0\text{ for any }a\ing.
\end{equation}

If $\f00=0$ then $\f a0=0=\f00$ by (\ref{eq3.04}). The lemma stands.

If $\f00\neq0$, then we only need to prove $\f a0\neq0$ by (\ref{eq3.04}).
Set
$$G_1=\{a\ing\mid\f a0=0\};\ \
  G_2=\{a\ing\mid\f a0\neq0\}=\{a\ing\mid\f a0=\f00\}.$$
Clearly, we have
$$G_1\cup G_2=G, G_1\cap G_2=\emptyset\text{ and }0\in G_2.$$
Set $c=0$ in (\ref{eq3.03}) and we have
\begin{equation}\label{eq3.05}
 (b-a)\f{a+b}0=\f b0\f ab-\f a0\f ba.
\end{equation}
Then we get the following conclusions
\begin{itemize}
 \item[(I)] If $a,b\in G_i, a\neq b,i=1,2$, then $a+b\in G_i$.
 \item[(II)] If $a\in G_1$, $-a\in G_2$.
\end{itemize}
Indeed, (I) follows directly from (\ref{eq3.05}) and the definition of $G_1,G_2$.
For (II) let $b=-a$ in (\ref{eq3.05}) and we get
$\f{-a}0\f a{-a}=-2a\f00\neq0$,
which implies that $\f{-a}0\neq 0$. So $-a\in G_2$.

{\bf Claim}: $G_2=G$, or equivalently, $G_1=\emptyset$.
Obviously the lemma follows from this claim.
Suppose otherwise that there exists some nonzero $\theta\in G_1$.
Hence $-\theta\in G_2$ by the conclusion (II).
Moreover, we have $2\theta\in G_1$.
Otherwise, $\theta=-\theta+2\theta\in G_2$ by conclusion (I), a contradiction.
Finally, by induction one can prove that $n\theta\in G_1$ for any integer $n>0$.
Denote by $H$ the cyclic group generated by $\theta$.
Thus we have
$$G_1\cap H=\{\theta,2\theta,3\theta,\cdots\}\text{ and }
  G_2\cap H=\{0,-\theta,-2\theta,-3\theta,\cdots\}.$$

Let $n$ be a negative integer.
Set $a=\theta, b=n\theta$ in (\ref{eq3.05}) and we get
$$\theta (n-1)\f00=\f00\f{\theta}{n\theta}.$$
Hence by (\ref{eq3.02})
\begin{equation}\label{eq3.06}
 \f{\theta}{n\theta}=(n-1)\theta \text{ and }\f{n\theta}{\theta}=0\text{ for any } n<0.
\end{equation}
On the other hand, choose $m,n\in\Z$ such that $m\geq-n>0$,
then $\f{n\theta}0=\f00\neq0$ and $\f{(m+1+n)\theta}0=0$.
Let $a=(m+1)\theta, b=n\theta$ in (\ref{eq3.05}) and we get
\begin{equation}\label{eq3.07}
 \f{(m+1)\theta}{n\theta}=0 \text{ and }\f{n\theta}{(m+1)\theta}=(m+n+1)\theta\text{ if } m\geq-n>0.
\end{equation}
Now let $a=m\theta,b=n\theta,c=\theta$ in (\ref{eq3.05}) and we have
$$(n-m)\theta\f{m\theta+n\theta}\theta=\f{n\theta}\theta\f{m\theta}{(n+1)\theta}-
  \f{m\theta}\theta\f{n\theta}{(m+1)\theta},$$
that is,
\begin{equation}\label{eq3.08}
 (n-m)\f{m\theta+n\theta}\theta=(n-m-1)\f{m\theta}\theta\text{ if }m\geq-n>0.
\end{equation}
Especially, let $m=-n>0$, then we have
\begin{equation}\label{eq3.09}
 \f{m\theta}\theta=\frac{2m}{2m+1}\f0\theta=\frac{2m}{2m+1}\theta\text{ for any }m\geq1.
\end{equation}
Here we have used $\f0\theta=\theta$.

Moreover, by taking $m=1-n>1$ in (\ref{eq3.08}) we get
$$(1-2m)\f\theta\theta=-2m\f{m\theta}\theta.$$
So by (\ref{eq3.09}) we obtain
$$\f\theta\theta=\frac{4m^2\theta}{4m^2-1}\text{ for any }m\geq2.$$
This is a contradiction since $\f\theta\theta$ is independent of $m$.
So the claim, hence the lemma, holds.
}

Since $L_0L_a=\f0aL_a=(\f00+a)L_a$,
the compatible left-symmetric algebra $W$ is a direct sum of $L_0$-eigenspaces
$$W=\bigoplus_{a\ing}W_a,$$
where $W_a=\{v\in W\mid L_0v=(\f00+a)v\}=\C L_a$.
Consider $W$ as a $\W$-module, this is also the weight space decomposition of $W$.

\begin{prop}\label{prop3.2}
As a $\W$-module, $W$ is indecomposable,
hence must be isomorphic to one of $V_{\al,\be},A_\gm,B_\gm$.
\end{prop}
\pf{
Assume $W=W_1\oplus W_2$, where $W_1,W_2$ are nonzero proper $\W$-submodules of $W$.
Without loss of generality, we may assume that
$W_2$ is indecomposable and $L_0\notin W_2$.

Since a submodule of a weight module is still a weight module
and all weight spaces of $W$ is 1-dimensional,
there exist nonempty subsets $S_1,S_2$ of $G$ such that
\begin{center}
 $0\notin S_2, S_1\cap S_2=\emptyset, S_1\cup S_2=G$
 and $W_i=\bigoplus_{a\in S_i}\C L_a,\ \ \ i=1,2$.
\end{center}
Then the following properties of $S_1, S_2$ and $G$ hold.
\begin{itemize}
 \item[(i)] $\f a0=0, \f0a=a$ for any $a\ing$.
 Take $b\in S_2$. Since $L_bL_0=\f b0L_b\in W_1\cap W_2=\{0\}$,
 we see $\f b0=0$. Hence $\f a0=\f00=\f b0=0$ by Lemma \ref{lem3.1}.

 \item[(ii)] If $a,b\in S_2$ and $a\neq b$, then $a+b\in S_2$.
 Since $\f ab-\f ba=b-a\neq 0$, at least one of $\f ab,\f ba$ is not zero.
 Let $\f ab\neq 0$. Then $L_aL_b=\f abL_{a+b}\neq 0$.
 So $L_{a+b}\in W_2, a+b\in S_2$.

 \item[(iii)] If $a\in S_2$ then $-a\in S_1$.
 Otherwise, $0=a+(-a)\in S_2$ by (ii), a contradiction.

 \item[(iv)] If $a\in S_i, \f ba\neq 0$ then $a+b\in S_i$ (i=1,2).
 This follows from $L_aL_b=\f abL_{a+b}$. So does the next property (v).

 \item[(v)] Let $a\in S_2, b\in S_1$. We have $a+b\in S_1$ if and only if $\f ba=0$,
 and $a+b\in S_2$ if and only if $\f ab=0$.
\end{itemize}

Since $W_2\neq \{0\}$, there exists $0\neq \theta\ing$ such that $\Z\theta\cap S_2\neq\emptyset$.\\
{\bf Case {I}}: $\Z_+\theta\cap S_2\neq\emptyset$.
Let $s$ be the least positive integer such that $s\theta\in S_2$.
Notice that $\f0{s\theta}=s\theta$
and $f(-\theta,s\theta)=0$ by property (iv) and the choice of $s$.
Let $a=-\theta, b=\theta, c=s\theta$ in (\ref{eq3.03}) and we get
\begin{equation}\label{eq3.10}
 0\neq 2s\theta^2=2\theta\f0{s\theta}
 =f(\theta,s\theta)\f{-\theta}{(s+1)\theta}-f(-\theta,s\theta)\f{\theta}{(s-1)\theta},
\end{equation}
from which we have $f(\theta,s\theta)\neq0$. So $(s+1)\theta\in S_2$ by (iv).

If $s=1$, then $2\theta\in S_2$. Moreover $\Z_+\theta\subseteq S_2$ by (ii),
and $\{0\}\cup\Z_-\theta\subseteq S_1$ by (iii).
Since $-\theta\in S_1$ and $\theta,2\theta\in S_2$, we have $f(2\theta,-\theta)=0$ by (v).
Hence $f(-\theta, 2\theta)=3\theta$.
Then from (\ref{eq3.10}) one gets $\f\theta\theta=\frac{2\theta}3$.
Similarly we have $f(2\theta,-\theta)=0, f(3\theta,-2\theta)=0$ and $f(-2\theta,3\theta)=5\theta$.
Let $a=-b=2\theta, c=\theta$ in (\ref{eq3.03}) and we see
$\f\theta\theta=\frac{4\theta}5$.
This forces $\theta=0$, a contradiction.

Suppose $s>1$. We have $(s-1)\theta,\theta\in S_1$.
Hence $L_\theta L_{(s-1)\theta}\in W_1\cap W_2=\{0\}$ and $L_{(s-1)\theta}L_\theta=0$.
This implies $\f\theta{(s-1)\theta}=\f{(s-1)\theta}\theta=0$.
On the other hand we have $\f\theta{(s-1)\theta}-\f{(s-1)\theta}\theta=(s-2)\theta$.
So $s=2$.

Notice that $-\theta\in S_1$ (otherwise, $\theta=-\theta+2\theta\in S_2$ by (ii), a contradiction),
$-2\theta\in S_1$ by (iii) and $3\theta\in S_2$.
Then by (v) we get $\f{2\theta}\theta=0$ and $\f\theta{2\theta}=\theta$.
Let $a=-b=2\theta, c=2\theta$ in (\ref{eq3.03}) and we have
$$-4\theta\f0{2\theta}=\f{-2\theta}{2\theta}\f{2\theta}0-\f{2\theta}{2\theta}\f{-2\theta}{4\theta},$$
which implies $\f{2\theta}{2\theta}\neq 0$.
Then $4\theta=2\theta+2\theta\in S_2$ by (iv).
Furthermore, we have $\f{4\theta}{-2\theta}=0$ by (v), hence $\f{-2\theta}{4\theta}=6\theta$.
So $\f{2\theta}{2\theta}=\frac{4\theta}3$.

Now let $a=2\theta, b=-\theta, c=2\theta$ in (\ref{eq3.03}) and
notice that $\f{-\theta}{4\theta}=5\theta$ by (v).
We get
$$-3\f\theta{2\theta}=\f{-\theta}{2\theta}\f{2\theta}\theta-\f{2\theta}{2\theta}\f{-\theta}{4\theta},$$
which implies $3=\frac{20}3$. This contradiction invalidates Case I.\\
{\bf Case II}: $\Z_-\theta\cap S_2\neq \emptyset$.
The discussion is similar to Case I and we omit it.
}

\begin{thm}\label{thm3.3}
 As a $\W$-module, $W$ is not isomorphic to $A_\gm$ for any $\gm\in\C$.
\end{thm}
\pf{
Suppose that $W\cong A_\gm$ for some $\gm\in\C$,
and $g:W\longrightarrow A_\gm$ is the isomorphism map.
Let $g(L_0)=\sum_{a\ing}C_av_a$.
Since $g(L_0L_0)=L_0g(L_0)$, we get
$$\sum_{a\ing}\f00C_av_a=\sum_{a\neq0}C_av_a.$$
If there exists some $a\neq 0$ such that $C_a\neq 0$,
then $\f00=a$ and $g(L_0)=C_av_a$.
By Lemma \ref{lem3.1} we know $\f b0=a$ for any $b\ing$.
Consider $g(L_{-a}L_0)=L_{-a}g(L_0)$ and we get
$$\f{-a}0g(L_{-a})=C_aL_{-a}v_a=0.$$
Therefore $g(L_{-a})=0$, which contradicts to that $g$ is an isomorphism.
So $g(L_0)=C_0v_0$, $C_0\neq 0$ and $\f00=0$.
Then for any $a\ing$ such that $a\neq0,-\gm$, we have
$$0=\f a0g(L_a)=g(\f a0L_a)=g(L_aL_0)=L_ag(L_0)=C_0L_av_0=C_0a(a+\gm)v_a.$$
This is also a contradiction.
}

\begin{thm}\label{thm3.4}
The following equation defines a graded CLSAS $V_{\al,\theta}$ on $\W$
\begin{equation}\label{eq3.11}
 L_bL_c=\frac{(\al+c+\al\theta b)(1+\theta c)}{1+\theta(b+c)}L_{b+c},
\end{equation}
where $\al,\theta\in\C$ satisfying that $\theta=0$ or $\theta^{-1}\not\ing$.
Furthermore, if a graded CLSAS $W$ on $\W$ is isomorphic to $V_{\al,\be}$ as a $\W$-module,
then it must be defined by (\ref{eq3.11}).
\end{thm}
\pf{
One can directly check that (\ref{eq3.11}) defines a graded CLSAS on $\W$.

Conversely, let $W\cong V_{\al,\be}$ as $\W$-modules.
Denote by $g:W\longrightarrow V_{\al,\be}$ the isomorphism map and write
\begin{equation}\label{eq3.12}
 g(L_0)=\sum_{a\ing}C_av_a,
\end{equation}
where the sum takes over a finite number of $a\ing$.

Notice that $V_{\al,\be}\cong V_{\al',\be'}$ as $\W$-modules if and only if
(i) $\al-\al'\ing,\be=\be'$, or (ii) $\al-\al'\ing,\al\not\ing$ and $\{\be,\be'\}=\{0,1\}$
(see \cite{Su} Proposition 2.2).
Recall the $\Z$-basis $\e_1,\cdots,\e_n$ of $G$.
In the following, we may assume that if $\al\in\R G$ and
write $\al=x_1\e_1+\cdots+x_n\e_n$ for $x_1,\cdots,x_n\in\R$,
then $0\leq x_i<1$ for all $1\leq i\leq n$.

Consider $g(L_0L_0)=L_0g(L_0)$ and we get
\begin{equation}\label{eq3.13}
 \sum_{a\ing}\f00C_av_a=\sum_{a\ing}C_aL_0v_a=\sum_{a\ing}C_a(\al+a)v_a.
\end{equation}
Therefore the sum in (\ref{eq3.12}) has only one summand.
Hence there exists some $a\ing$ such that $g(L_0)=C_aL_a$.
So $\f00=\al+a$ by (\ref{eq3.13}) and $\f b0=\al+a$ for any $b\ing$ by Lemma \ref{lem3.1}.
Moreover from $g(L_bL_0)=L_gg(L_0)$ we obtain
\begin{equation}\label{eq3.14}
 \f b0g(L_b)=C_aL_bv_a=C_a(\al+a+b\be)v_{a+b}.
\end{equation}

If $\al+a\neq 0$ then
$$g(L_b)=\frac{\al+a+b\be}{\al+a}C_av_{a+b}.$$
Notice that $\al+a+b\be\neq 0$ for any $b\ing$ since $g$ is an isomorphism.
Then from $g(L_bL_c)=L_bg(L_c)$ we get
$$\f bc\frac{\al+a+(b+c)\be}{\al+a}C_av_{a+b+c}=\f bcg(L_{b+c})
  =\frac{\al+a+c\be}{\al+a}(\al+a+c+b\be)C_av_{a+b+c},$$
which gives
$$\f bc=\frac{(\al+a+c\be)(\al+a+c+b\be)}{\al+a+(b+c)\be}.$$
Replacing $\al+a$ and $\frac{\be}{\al+a}$ by $\al$ and $\theta$ respectively,
we get
\begin{equation}\label{eq3.15}
  \f bc=\frac{(\al+c+\al\theta b)(1+\theta c)}{1+\theta(b+c)}\ \ \ \ \
   \text{ if }\al+a\neq0.
\end{equation}
Notice that $g(L_b)=(1+\theta b)C_av_{a+b}$ by (\ref{eq3.14}).
Since $g$ is an isomorphism, it forces $\theta=0$ or $\theta^{-1}\not\ing$.

If $\al+a=0$, then for any $b\ing$ we have $\f b0=\f00=\al+a=0$.
Therefore (\ref{eq3.14}) forces $\be=0$.
Moreover we have $\al=-a\ing$ and by our assumption on $\al$ we see $\al=0=a$.

Since $L_0g(L_b)=g(L_0L_b)=\f 0bg(L_b)=bg(L_b)$, we see that
$$g(L_b)=A_bv_b\text{ for some }A_b\neq 0.$$
From $\f bcA_{b+c}v_{b+c}=\f bcg(L_{b+c})=g(L_bL_c)=L_bg(L_c)=A_ccv_{b+c}$ we get
\begin{equation}\label{eq3.16}
 \f bc=\frac{cA_c}{A_{b+c}},
\end{equation}
and moreover by (\ref{eq3.02})
\begin{equation}\label{eq3.17}
 cA_c-bA_b=(c-b){A_{b+c}}.
\end{equation}
In particular we have
\begin{equation}\label{eq3.18}
 A_b+A_{-b}=2A_0\text{ for any }b\ing.
\end{equation}
Replace $b$ by $-b$, $c$ by $2b$ in (\ref{eq3.17}) and apply (\ref{eq3.18}),
then we get $A_{2b}=2A_b-A_0$.
Moreover an induction shows that
\begin{equation}\label{eq3.19}
 A_{kb}=kA_b-(k-1)A_0\text{ for any }k\in\Z.
\end{equation}
For $b\ing$ such that $b\notin\Z\e_1$, from (using (\ref{eq3.17}) and (\ref{eq3.19}))
$$\begin{aligned}
 bA_b-\e_1A_{e_1}&=(b-\e_1)A_{b+\e_1}
   =\frac{b-\e_1}{3\e_1-b}\left(2\e_1A_{2\e_1}-(b-\e_1)A_{b-\e_1}\right)\\
   &=\frac{b-\e_1}{3\e_1-b}\left(2\e_1A_{2\e_1}-\frac{b-\e_1}{b+\e_1}(bA_b+\e_1A_{-\e_1})\right)\\
   &=\frac{b-\e_1}{3\e_1-b}\left(2\e_1(2A_{\e_1}-A_0)
          -\frac{b-\e_1}{b+\e_1}(bA_b-\e_1A_{\e_1}+2\e_1A_0)\right),
\end{aligned}$$
we obtain
$$A_b=\frac b{\e_1}A_{\e_1}-(\frac b{\e_1}-1)A_0=\frac b{\e_1}(A_{\e_1}-A_0)+A_0
 \ \ \ \ \text{ for any }b\ing.$$
Set $\theta=\frac{A_{\e_1}-A_0}{A_0\e_1}$ and we have by (\ref{eq3.16})
$$\f bc=\frac{cA_c}{A_{b+c}}=\frac{1+\theta c}{1+\theta(b+c)},$$
which coincides with (\ref{eq3.15}) with $\al=0$.
By a same argument as in the $\al+a\neq0$ case we also have $\theta=0$ or $\theta^{-1}\not\ing$.
}

\begin{thm}\label{thm3.5}
The following equation defines a graded CLSAS $V^{\gm,\lmd}$ on $\W$
 \begin{equation}\label{eq3.20}
  L_bL_c=\begin{cases}
                (\lmd+c)L_{b+c}  &\text{if }\lmd+b+c\neq0;\\
                \frac{(\lmd+c)(\gm-\lmd-c)}{\gm-\lmd}L_{-\lmd}&\text{if }\lmd+b+c=0,
         \end{cases}
 \end{equation}
 where $\gm\in\C,\lmd\ing$ such that $\gm-\lmd\neq 0$.
Furthermore, if a graded CLSAS $W$ on $\W$ is isomorphic to $B_\gm$ as a $\W$-module,
then it must be $V^{\gm,\lmd}$ defined by (\ref{eq3.20})
for some $\lmd\neq\gm$.
\end{thm}
\pf{
One can easily check that (\ref{eq3.20}) defines a graded CLSAS on $\W$.

Conversely, suppose that $W\cong B_\gm$ and
denote by $g:W\longrightarrow B_\gm$ be the isomorphism map with
$g(L_0)=\sum_{a\ing}C_av_a$.
Since $g(L_0L_0)=L_g(L_0)$ and $L_0v_0=0$, we have
$$\sum_{a\ing}\f00C_av_a=\sum_{a\neq0}C_aav_a.$$

{\bf Case I}: $C_0=0$.
There exists some $0\neq \lmd\ing$ such that $\f00=\lmd$ and $g(L_0)=C_\lmd v_\lmd$.
By lemma \ref{lem3.1} we have $\f b0=\lmd, \f 0b=\lmd+b$ for any $b\ing$.
From $\f b0g(L_b)=g(L_bL_0)=L_bg(L_0)=C_\lmd L_bv_\lmd$, one can get
$$g(L_b)=\begin{cases}
           C_\lmd v_{\lmd+b} &\text{if }\lmd+b\neq0;\\
           C_\lmd(\gm-\lmd)v_0&\text{if }\lmd+b=0.
  \end{cases}$$
Moreover, from $g(L_bL_c)=L_bg(L_c)$ we get
$$\f bcg(L_{b+c})=\begin{cases}
           C_\lmd L_bv_{\lmd+c} &\text{if }\lmd+c\neq0;\\
           C_\lmd(\gm-\lmd)L_bv_0&\text{if }\lmd+c=0.
  \end{cases}$$
Notice that $\gm-\lmd\neq0$ since $g$ is an isomorphism.

{\bf Subcase I.1}: $C_0=0$ and $\lmd+b+c\neq 0$.
If $\lmd+c\neq0$ then
$$\f bcC_\lmd v_{\lmd+b+c}=C_\lmd(\lmd+c)v_{\lmd+b+c}.$$
Thus $\f bc=\lmd+c$.
If $\lmd+c=0$, then $\f bc=\f b{-\lmd}=0$ since $g$ is an isomorphism.
So
\begin{equation}\label{eq3.21}
 \f bc=\lmd+c\ \text{ for all } b,c\ing\text{ such that } \lmd+b+c\neq 0.
\end{equation}

{\bf Subcase I.2}: $C_0=0$ and $\lmd+b+c=0$.
If $\lmd+c\neq0$ then
$$\f bcC_\lmd(\gm-\lmd)v_{\lmd+b+c}=C_\lmd L_bv_{\lmd+c}=C_\lmd(\lmd+c)(\gm-\lmd-c)v_0.$$
So $\f bc=\f{-\lmd-c}c=\frac{(\lmd+c)(\gm-\lmd-c)}{\gm-\lmd}$.
If $\lmd+c=0$ then $b=0$ and $\f bc=0$. So
\begin{equation}\label{eq3.22}
 \f bc=\frac{(\lmd+c)(\gm-\lmd-c)}{\gm-\lmd}
 \ \text{ for all } b,c\ing\text{ such that } \lmd+b+c= 0.
\end{equation}

{\bf Case II}: $C_0\neq0$.
Then we must have $\f00=0, g(L_0)=C_0v_0$
and $\f b0=0,\f 0b=b$ for any $b\ing$.
From $L_0g(L_b)=g(L_0L_b)=\f0bg(L_b)=bg(L_b)$ we see that
$$g(L_b)=K_bv_b\text{ where }K_b\neq 0\text { and }K_0=C_0\neq 0.$$
From $(c-b)K_{b+c}v_{b+c}=g([L_b,L_c])=L_bg(L_c)-L_cg(L_b)$ we have
\begin{align}
 &(c-b)K_{b+c}=cK_c-bK_b\ \ \ \ \text{ if }b+c\neq 0,\label{eq3.23}\\
 &(\gm-c)K_c+(\gm+c)K_{-c}=2K_0. \label{eq3.24}
\end{align}
Fix any $0\neq x\ing$. Let $b=2x,c=-x$ in (\ref{eq3.23}) we get
$$K_{-x}=3K_x-2K_{2x}.$$
Let $b=-2x,c=x$ in (\ref{eq3.23}) we get
$$K_{-2x}=4K_x-3K_{2x}.$$
Put these two equations into (\ref{eq3.24}) with $c=x$ and $c=2x$ separatively,
and we have
$$\begin{cases}
    &(\gm-x)K_x+(\gm+x)(3K_x-2K_{2x})=2K_0,\\
    &(\gm-2x)K_{2x}+(\gm+2x)(4K_x-3K_{2x})=2K_0,
  \end{cases}$$
from which we get $K_x=K_{2x}$ and $K_0=\gm K_x$.
This implies that $\gm\neq0$ and $K_b=\frac{K_0}\gm$ for any $b\ing\setminus\{0\}$.

Now from
$$\f bcK_{b+c}v_{b+c}=g(L_bL_c)=L_bg(L_c)=\begin{cases}
    K_ccv_{b+c}     &\text{ if }b+c\neq 0;\\
    K_cc(\gm-c)v_0  &\text{ if }b+c= 0,
   \end{cases}$$
we get that
$$\f bc=\begin{cases}
    c    &\text{ if }b+c\neq 0;\\
    \frac{c(\gm-c)}\gm  &\text{ if }b+c= 0,
   \end{cases}$$
which coincides with (\ref{eq3.21}) and (\ref{eq3.22}) with $\lmd=0$.
}

\begin{lem}\label{lem3.6}
Let $(W_1,\cdot), (W_2,\ast)$ be two graded CLSAS's on $\W$ defined respectively by
$$L_a\cdot L_b=f_1(a,b)L_{a+b}\text{ and }L_a\ast L_b=f_2(a,b)L_{a+b}.$$
If $(W_1,\cdot)$ is isomorphic to $(W_2,\ast)$, then
$$f_1(a,b)=f_2(a,b)\text{ or }f_1(a,b)=-f_2(-a,-b).$$
\end{lem}
\pf{
Let $\sgm:W_1\longrightarrow W_2$ be an isomorphism.
Clearly $\sgm$ is also a Lie algebra automorphism of $\W$.
Then it is easy to see that $\sgm(L_0)=\pm L_0$.

If $\sgm(L_0)=L_0$, then for any $a\ing$ we have $\sgm(L_a)=K_aL_a$ with some $K_a\neq 0$.
Expanding $\sgm([L_a,L_b])$ we get
$$K_aK_b=K_{a+b}.$$
Since $f_1(a,b)K_{a+b}L_{a+b}=\sgm(L_a\cdot L_b)=\sgm(L_a)\ast\sgm(L_b)=K_aK_bf_2(a,b)L_{a+b}$,
we see $f_1(a,b)=f_2(a,b)$.

If $\sgm(L_0)=-L_0$, then we may write $\sgm(L_a)=C_aL_a$ with some $C_a\neq 0$.
Expanding $\sgm([L_a,L_b])$ we get
$$C_aC_b=-C_{a+b}.$$
Since $f_1(a,b)C_{a+b}L_{a+b}=\sgm(L_a\cdot L_b)=\sgm(L_a)\ast\sgm(L_b)=C_aC_bf_2(-a,-b)L_{-(a+b)}$,
we see $f_1(a,b)=-f_2(-a,-b)$.
}

\begin{thm}\label{thm3.7}
Let $W_1,W_2$ be two graded CLSAS's on $\W$.
Then $W_1\cong W_2$ if and only if one of the following four cases stands
\begin{itemize}
 \item[(1)] $W_1=W_2$;
 \item[(2)] $\{W_1,W_2\}=\{V_{\al,0},V_{\al,\frac1\al}\}\text{ with } \al\not\ing$;
 \item[(3)] $\{W_1,W_2\}=\{V_{\al,\theta},V_{-\al,-\theta}\}\text{ with } \al\in\C,
               \theta=0\text{ or }\theta^{-1}\not\ing$;
 \item[(4)] $\{W_1,W_2\}=\{V^{\gm,0},V^{-\gm,0}\}\text{ with }\gm\neq0$.
\end{itemize}
\end{thm}
\pf{
Let $W_1\neq W_2$ and we continue to use the notations in Lemma \ref{lem3.6}.
By Proposition \ref{prop3.2}, Theorem \ref{thm3.3}, \ref{thm3.4} and \ref{thm3.5},
we know $W_1,W_2$ must be $V_{\al,\theta}$ or $V^{\gm,\lmd}$.
Hence there are three cases to be considered.

{\bf Case I}: $W_1=V_{\al_1,\theta_1}, W_2=V_{\al_2,\theta_2}$.
If $f_1(b,c)=f_2(b,c)$, then $\al_1=f_1(0,0)=f_2(0,0)=\al_2$.
Write $\al=\al_1=\al_2$.
Notice that $\theta_1\neq\theta_2$ since $W_1\neq W_2$.
From
$$f_1(b,c)=\frac{(\al+c+\al\theta_1b)(1+\theta_1c)}{1+\theta_1(b+c)}
  =f_2(b,c)=\frac{(\al+c+\al\theta_2b)(1+\theta_2c)}{1+\theta_2(b+c)}$$
we see that
$$bc(\theta_1-\theta_2)\left(\al(\theta_1+\theta_2)-1+\al\theta_1\theta_2(b+c)\right)=0,$$
which implies that $\al\neq 0, \theta_1\theta_2=0$ and $\theta_1+\theta_2=\frac1\al$.
This is case (2).
We may assume $\theta_1=0$ by symmetry, then
we have $\al=\theta_2^{-1}\not\ing$ in case (2).

If $f_1(b,c)=-f_2(-b,-c)$, then we have $\al_1=-\al_2$.
Write $\al=\al_1=-\al_2$. From
$$f_1(b,c)=\frac{(\al+c+\al\theta_1b)(1+\theta_1c)}{1+\theta_1(b+c)}
  =f_2(b,c)=\frac{(\al+c-\al\theta_2b)(1-\theta_2c)}{1-\theta_2(b+c)}$$
we see that
$$bc(\theta_1+\theta_2)\left(\al(\theta_1-\theta_2)-1-\al\theta_1\theta_2(b+c)\right)=0,$$
which gives $\theta_1=-\theta_2$ or
$\theta_1\neq-\theta_2, \theta_1\theta_2=0$ and $\theta_1-\theta_2=\frac1\al$.
Hence
\begin{align*}
 &W_1=V_{\al,\theta},\ W_2=V_{-\al,-\theta}, \text{ with } \theta=0\text{ or }\theta^{-1}\not\ing;\\
 \text{or }&W_1=V_{\al,0},\ W_2=V_{-\al,-\frac1\al}\text{ with } \al\not\ing;\\
 \text{or }&W_1=V_{\al,\frac1\al},\ W_2=V_{-\al,0}\text{ with } \al\not\ing.
\end{align*}
Applying case (2) we can see that what the above equations describe is exactly case (3).

{\bf Case II}: $W_1=V^{\gm_1,\lmd_1},\ W_2=V^{\gm_2,\lmd_2}$.
If $f_1(b,c)=f_2(b,c)$, then $\lmd_1=f_1(0,0)=f_2(0,0)=\lmd_2$.
Take $a\neq 0,-\lmd_1$. Since
$$f_1(-\lmd_1-a,a)=\frac{(\lmd_1+a)(\gm_1-\lmd_1-a)}{\gm_1-\lmd_1}
  =f_2(-\lmd_2-a,a)=\frac{(\lmd_2+a)(\gm_2-\lmd_2-a)}{\gm_2-\lmd_2},$$
we get $\gm_1=\gm_2$, contradicting to $W_1\neq W_2$.
Hence $f_1(b,c)=-f_2(-b,-c)$ and we have $\lmd_1=-\lmd_2$.
Write $\lmd=\lmd_1=-\lmd_2$.
Suppose $\lmd\neq0$. Choose $a\neq 0, -\lmd$ and then from
$f_1(-\lmd_1-a,a)=-f_2(\lmd_2+a,-a)$ we get
$$\frac{(\lmd+a)(\gm_1-\lmd-a)}{\gm_1-\lmd}=\lmd+a,$$
which is a contradiction. Therefore $\lmd=0$.
Now from
$$f_1(-a,a)=\frac{a(\gm_1-a)}{\gm_1}=-f_2(a,-a)=-\frac{-a(\gm_2+a)}{\gm_2}$$
it follows that $\gm_1=-\gm_2$. This is the case (4).

{\bf Case III}. $W_1=V_{\al,\theta},\ W_2=V^{\gm,\lmd}$.
By case (2) we may replace $W_1$ by $V_{-\al,-\theta}$
and suppose that $f_1(b,c)=-f_1(-b,-c)$.
Therefore we always have $f_1(b,c)=f_2(b,c)$.
Thus $\al=f_1(0,0)=f_2(0,0)=\lmd\ing$.

Now choose $b,c\neq0$ such that $\lmd+b+c\neq0$. Since
$$f_1(b,c)=\frac{(\lmd+c+\theta\lmd b)(1+\theta c)}{1+\theta(b+c)}
  =f_2(b,c)=\lmd+c,$$
we get $\theta bc(\lmd\theta-1)=0$.
Then we must have $\theta=0$.
Otherwise we have $\theta^{-1}=\lmd\ing$, a contradiction.

Now let $c\neq 0,-\lmd$. Then we have
$$0\neq f_1(-\lmd-c,c)=\lmd+c=f_2(-\lmd-c,c)=\frac{(\lmd+c)(\gm-\lmd-c)}{\gm-\lmd},$$
which is also a contradiction.
So there is no isomorphism between $V_{\al,\theta}$ and $V^{\gm,\lmd}$.
}

\section{Non-graded CLSAS's on $\W$}\label{sec4}
\def\theequation{4.\arabic{equation}}
\setcounter{equation}{0}

In this section we consider CLSAS's on $\W$ with multiplication
\begin{equation}\label{eq4.01}
 L_aL_b=\f abL_{a+b}+\g abL_{a+b+\zeta},
\end{equation}
where $\zeta\ing$ is a fixed nonzero number,
and $f,g$ are two complex-valued functions on $G\times G$.

Clearly if $g$ is not a zero map,
then the CLSAS on $\W$ satisfying (\ref{eq4.01}) is non-graded.
Let $\al,\mu\in\C$ and set
$$\f ab=\al+b, \g ab=\mu\text{ for any }a,b\ing.$$
One can easily check that this defines a CLSAS on $\W$
and it is non-graded if $\mu\neq0$.
We denote such a CLSAS by $W(\al,\mu,\zeta)$.
Moreover, we have $V_{\al,0}=W(\al,0,\zeta)$.
Our main result in this section is the following

\begin{thm}\label{thm4.1}
Any non-graded CLSAS on $\W$ defined by (\ref{eq4.01}) is isomorphic to some $W(\al,\mu,\zeta)$
where $\al\in\C,\mu\in\C\setminus\{0\}$.
\end{thm}

Before giving the proof of the theorem, we should state some properties about $f,g$.
\begin{lem}\label{lem4.2}
The equation (\ref{eq4.01}) defines a CLSAS on $\W$ if and only if
$f,g$ satisfy (\ref{eq3.02}), (\ref{eq3.03}) and the following three equations
\begin{align}
 &\g ab=\g ba; \label{eq4.02}\\
 &\g bc\g a{b+c+\zeta}=\g ac\g b{a+c+\zeta}; \label{eq4.03}\\
 &\f bc\g a{b+c}-\f ac g(b,a+c)\notag \\
 &=\g ac f(b,a+c+\zeta)-\g bc f(a,b+c+\zeta)+(b-a)g(a+b,c).\label{eq4.04}
\end{align}
\end{lem}
\pf{
This lemma follows directly from
$$\begin{aligned}
  &[L_a,L_b]=L_aL_b-L_bL_a=(b-a)L_{a+b}\\
  &(L_aL_b)L_c-(L_bL_a)L_c=L_a(L_bL_c)-L_b(L_aL_c).
 \end{aligned}$$
}

Since the equations (\ref{eq3.02}) and (\ref{eq3.03}) define a graded CLSAS on $\W$,
by Theorem \ref{thm3.4} and Theorem \ref{thm3.5} we know that
\begin{equation}\label{eq4.05}
 \f ab=\frac{(\al+b+a\al\theta)(1+\theta b)}{1+\theta(a+b)}
\end{equation}
with $\al,\theta\in\C$ such that $\theta=0$ or $\theta^{-1}\not\ing$, or
\begin{equation}\label{eq4.06}
 \f ab=\begin{cases}
  \lmd+b &\text{ if } \lmd+a+b\neq0;\\
  \frac{(\lmd+b)(\gm-\lmd-b)}{\gm-\lmd} &\text{ if } \lmd+a+b=0
  \end{cases}
  =(\lmd+b)\left(1-\frac b{\gm-\lmd}\dt_{\lmd+a+b,0}\right)
\end{equation}
with $\gm\in\C,\lmd\ing$ such that $\gm\neq\lmd$.

\begin{prop}\label{prop4.3}
If the map $f$ takes the form (\ref{eq4.06}), then $g=0$.
\end{prop}
\pf{
Note that $\f a0=\lmd,\f0a=\lmd+a$ for any $a\ing$.
Take $a=c=0$ in (\ref{eq4.04}) and we get
$$(\lmd+\zeta)\g b0=(\lmd+\zeta)\left(1-\frac\zeta{\gm-\lmd}\dt_{\lmd+\zeta+b,0}\right)\g00.$$

{\bf Case I}: $\lmd+\zeta\neq0$.
The above equation turns to
\begin{equation}\label{eq4.07}
 \g b0=\left(1-\frac\zeta{\gm-\lmd}\dt_{\lmd+\zeta+b,0}\right)\g00.
\end{equation}
Let $a=0$ in (\ref{eq4.04}) and we have
\begin{equation}\label{eq4.08}
 \zeta\g bc=\f b{c+\zeta}\g c0-\f bc\g{b+c}0.
\end{equation}
Notice by (\ref{eq4.06}), (\ref{eq4.07}) and (\ref{eq4.08}) that
$$g(2\zeta,-\lmd+\zeta)=g(2\zeta,-\lmd-\zeta)=g(\zeta,2\zeta-\lmd)=\g00,\ \
  g(\zeta,-\lmd-\zeta)=\frac{\gm+\zeta}{\gm-\zeta}\g00.$$
Now set $a=2\zeta, b=\zeta, c=-\lmd-\zeta$ in (\ref{eq4.03}) and we obtain
$$\frac{\gm+\zeta}{\gm-\zeta}g^2(0,0)=g^2(0,0),$$
which forces $\g00=0$ since $\lmd+\zeta\neq0$.
So it follows from (\ref{eq4.07}) and (\ref{eq4.08}) that
$\g bc=0$ for all $b,c\ing$.

{\bf Case II}: $\lmd+\zeta=0$.
Let $a=0,c=\zeta$ in (\ref{eq4.04}) and we get
\begin{equation}\label{eq4.09}
\g b\zeta=\left(1-\frac{2\zeta}{\gm+\zeta}\dt_{b+\zeta,0}\right)\g\zeta0.
\end{equation}
Take $b=\zeta,c=0$ and $a\neq\zeta$ in (\ref{eq4.04}), then we have
\begin{equation}\label{eq4.10}
 a\g a0+(\zeta-a)g(a+\zeta,0)=\zeta\g\zeta0\ \text{ if }a\neq-\zeta.
\end{equation}
Taking $b=c=0$ in (\ref{eq4.03}), it gives
\begin{equation}\label{eq4.11}
 \g a0g(a+\zeta,0)=\g00\g a\zeta=\g00\g\zeta0\ \text{ if }a\neq-\zeta.
\end{equation}
Combing (\ref{eq4.10}) and (\ref{eq4.11}), we have
\begin{equation}\label{eq4.12}
 ag^2(a,0)-\zeta g(\zeta,0)\g a0+(\zeta-a)\g00\g\zeta0=0\ \text{ if }a\neq-\zeta.
\end{equation}
Replace $a$ by $a+\zeta$ in (\ref{eq4.12}) and we get
\begin{equation}\label{eq4.13}
 (a+\zeta)g^2(a+\zeta,0)-\zeta g(\zeta,0)\g{a+\zeta}0-a\g00\g\zeta0=0\ \text{ if }a\neq-2\zeta.
\end{equation}

{\bf Claim}: $\g\zeta0=0$.\\
Suppose $\g\zeta0\neq0$.
We may assume $\g\zeta0=1$ by replacing $\g ab$ with $\g ab/\g\zeta0$
in equations (\ref{eq4.02})-(\ref{eq4.04}), if necessary.
Let $a\neq\zeta$,
multiply $(\zeta-a)^2$ to (\ref{eq4.13}) and combine (\ref{eq4.02}) and (\ref{eq4.03}),
then one can get
\begin{equation}\label{eq4.14}
 a^2 g(a,0)+a(\zeta-a)\g00-a\zeta=0\ \text{ if }a\neq\pm\zeta,-2\zeta.
\end{equation}
Moreover, since (here we use equations (\ref{eq4.10})-(\ref{eq4.13}))
$$\begin{aligned}
 &a(a+\zeta) g^2(0,0)=a(a+\zeta) g^2(a,0) g^2(a+\zeta,0)\\
 &=\left(\zeta\g a0-(\zeta-a)\g00\right)\left(\zeta g(a+\zeta,0)+a\g00\right)\\
 &=\zeta^2\g00+a\zeta \g00\g a0-\zeta \g00(\zeta-a\g a0)-a(\zeta-a)g^2(0,0)\\
 &=2a\zeta \g00\g a0-a(\zeta-a)g^2(0,0),
 \end{aligned}$$
we see that
\begin{equation}\label{eq4.15}
 g^2(0,0)=\g00\g a0\text{ if }a\neq 0,-\zeta,-2\zeta.
\end{equation}
Put equation (\ref{eq4.14}) into (\ref{eq4.15}), then we get
$$\g00(\g00-1)=0.$$
If $\g00=0$, then by (\ref{eq4.14}) we have
$$\g a0=\frac\zeta a\text{ for all }a\neq0,\pm\zeta,-2\zeta.$$
Let $a=2\zeta$ in (\ref{eq4.11}) and we get
$$0=g(2\zeta,0)g(3\zeta,0)=\frac16.$$
This contradiction forces $\g00=1$.
Then by (\ref{eq4.15}) $\g a0=1$ for any $a\neq-\zeta,-2\zeta$.
Take $a=-2\zeta$ and $a=-3\zeta$ in (\ref{eq4.10}) respectively and we obtain
$g(-\zeta,0)=g(-2\zeta,0)=1$.
Hence $\g a0=1$ for all $a\ing$.

Now let $a=0$ in (\ref{eq4.04}) and we have
\begin{equation}\label{eq4.16}
 \zeta\g bc=\f b{c+\zeta}-\f bc,
\end{equation}
which implies $\g bc=1$ for all $b,c$ such that $b+c\neq-\lmd,-\lmd-\zeta$.
If $b+c=-\lmd$ or $-\lmd-\zeta$,
take $a\ing$ such that $a+c\neq0,\zeta$ and $a+b+c\neq0,-\zeta$ in (\ref{eq4.03}),
then we have $\g bc=1$ for all $b,c\ing$.
Therefore (\ref{eq4.16}) becomes
$$\f b{c+\zeta}-\f bc=\zeta.$$
But when take $b,c\neq0$ such that $\lmd+b+c=0$ in (\ref{eq4.06}), one gets
$$\f b{c+\zeta}-\f bc=\zeta-\frac{bc}{\gm-\lmd}.$$
With this contradiction we prove the claim that $\g\zeta0=0$.

By (\ref{eq4.09}) we have $\g b\zeta=0$ for all $b\ing$,
and then from (\ref{eq4.12}) we see that
$$\g a0=0\text{ for }a\neq 0,-\zeta.$$
Take $a=-2\zeta$ in (\ref{eq4.10}) and we get $\g{-\zeta}0=0$.
Let $a=-b\neq0, c=0$ in (\ref{eq4.04}) and we see $\g00=0$.
So $\g a0=0$ for any $a\ing$.
Then take $a=0$ in (\ref{eq4.04}) and we have
$\g bc=0$ for all $b,c\ing$.
}

Before dealing with the case when $f$ takes the form (\ref{eq4.05}),
we shall recall some results about symmetric rational functions from \cite{TB}.
\begin{prop}[\cite{TB} Proposition 2.6]\label{prop4.4}
Let $h(a,b)$ be a symmetric rational function in variables $a,b\in\C$.
Then $h(a,b)$ satisfies equation (\ref{eq4.03}) if and only if
$h(a,b)=\frac{T(a)T(b)}{T(a+b+\zeta)}$ for some rational function $T(x)$.\\
\end{prop}
\begin{prop}[\cite{TB} Proposition 2.7]\label{prop4.5}
Suppose $T(x)$ is a nonzero symmetric rational function on $\C$ and $\eta\in\C\setminus\{0\}$.
\begin{itemize}
 \item[(1)] If $\frac{T(x)}{T(x+\eta)}=\frac{\sum_{i=0}^sa_ix^i}{\sum_{i=0}^tb_ix^i}$
      with $a_sb_t\neq0$, then $s=t$ and $a_s=b_t$.
      In particular, $\frac{T(x)}{T(x+\eta)}\neq(a+bx)^{\pm1}$ for any $a,b\in\C$ with $b\neq0$.

 \item[(2)] If $\frac{T(x)}{T(x+\eta)}=\chi$ is a constant,
       then $\chi=1$ and $T(x)$ is also a constant.
 \item[(3)] If $\frac{T(x)}{T(x+\eta)}=\rho\frac{x+a}{x+b}$ with $a,b,\rho\in\C$ such that
       $a\neq b$ and $\rho\neq0$, then $\rho=1$ and there exist some $k\in\Z\setminus\{0\}$
       and $\tau\in\C\setminus\{0\}$ such that $b-a=k\eta$ and
       $$ T(x)=\begin{cases}
             \tau\prod\limits_{j=0}^{k-1}(x+a+j\eta)&\text{ if }k\geq1;\\
             \tau\prod\limits_{j=0}^{-k-1}(x+b+j\eta)^{-1}&\text{ if }k\leq-1.
           \end{cases}$$
\end{itemize}
\end{prop}

We note that Proposition \ref{prop4.5} is originally presented in \cite{TB} with $\eta\in\Z$.
But the proof still goes through if we replace $\eta>0$ with $\text{Re}\eta>0$ for $\eta\in\C$.

\begin{prop}\label{prop4.6}
If the map $f$ takes the form (\ref{eq4.05}) and $g\neq0$,
then $\theta=0, f(a,b)=\al+b$ and $g$ is a constant.
\end{prop}
\pf{
Take $b=c=0$ in (\ref{eq4.04}) and use (\ref{eq4.02}), (\ref{eq4.05}), then we get
\begin{equation}\label{eq4.17}
 (\al+\zeta)\g a0=\frac{(\al+\zeta+a\al\theta)(1+\theta\zeta)}{1+\theta(a+\zeta)}\g00.
\end{equation}

{\bf Case I}: $\al+\zeta=0$.
Suppose $\theta\neq0$.
Then (\ref{eq4.17}) implies that $\g00=0$.
Let $b=c=0$ in (\ref{eq4.03}), we have
\begin{equation}\label{eq4.18}
 g(a+\zeta,0)\g a0=0\text{ for any }a\ing.
\end{equation}
Set $a=-\zeta, c=0$ in (\ref{eq4.04}) and we obtain
\begin{equation}\label{eq4.19}
 \zeta g(-\zeta,0)+\frac{(b+\theta\zeta^2)(1+\theta(b+\zeta))}{1+\theta b}\g b0
 =(b+\zeta)g(b-\zeta,0).
\end{equation}
Put $b=\zeta$ in (\ref{eq4.19}), we get
\begin{equation}\label{eq4.20}
 g(-\zeta,0)=-(1+2\theta\zeta)\g\zeta0.
\end{equation}
Put $b=2\zeta$ in (\ref{eq4.19}) and use (\ref{eq4.20}), then we get
$$\frac{(2+\theta\zeta)(1+3\theta\zeta)}{1+2\theta\zeta} g(2\zeta,0)=2(2+\theta\zeta)\g\zeta0.$$
If $\theta\zeta\neq-2$, then $g(2\zeta,0)=\frac{2(1+2\theta\zeta)}{1+3\theta\zeta}\g\zeta0$.
Put into (\ref{eq4.18}) with $a=\zeta$ and we get $\g\zeta0=0$.
If $\theta\zeta=-2$, put $b=3\zeta$ in (\ref{eq4.19}), then we get
$$\frac75 g(3\zeta,0)=4g(2\zeta,0)-3\g\zeta0.$$
Multiply $g(2\zeta,0)$ and we see that $g(2\zeta,0)=0$.
Hence $g(3\zeta,0)=-\frac{15}7g(\zeta,0)$.
Now put $b=4\zeta$ in (\ref{eq4.19}) and we see $g(4\zeta,0)=-\frac{111}{18}\g\zeta0$.
Then by (\ref{eq4.18}) we still have $\g\zeta0=0$.

By (\ref{eq4.20}) we see $g(-\zeta,0)=0$.
Multiply $\g b0$ to (\ref{eq4.19}) and use (\ref{eq4.18}), then we get
$$\g b0=0\text{ for all }b\neq-\theta\zeta^2.$$
If $b=-\theta\zeta^2$, put $a=\theta\zeta^2, b=-2\theta\zeta^2, c=0$ in (\ref{eq4.04}),
we see $g(-\theta\zeta^2,0)=0$.
Thus $\g b0=0$ for all $b\ing$.
Set $a=0$ in (\ref{eq4.04}) and we see that
$\g bc=0$ for all $b,c\ing$,
contradicting to $g\neq0$.
So $\theta=0$ and $f(a,b)=\al+b$.

Now let $c=0$ in (\ref{eq4.04}) and we get
\begin{equation}\label{eq4.21}
 a\g a0-b\g b0+(b-a)g(a+b,0)=0.
\end{equation}
Furthermore, put $a+b=0$ and we see
\begin{equation}\label{eq4.22}
\g a0+g(-a,0)=2\g00\ \ \ \text{ for all }a\ing.
\end{equation}
Fix any $\xi\ing\setminus\{0\}$.
Let $a=2\xi,b=-\xi$ in (\ref{eq4.21}) and by (\ref{eq4.22}) we obtain
$$2\g\xi0-g(2\xi,0)=\g00.$$
Put $a=\xi,b=k\xi,k\in\Z_+$ in (\ref{eq4.21}), we get
$$(k-1)g((k+1)\xi,0)-kg(k\xi,0)+\g\xi0=0,$$
whose solution is $g(k\xi,0)=\tau k+\mu,\ k\in\Z_+$,
where $\tau=g(2\xi,0)-g(\xi,0)$ and $\mu=2\g\xi0-g(2\xi,0)=\g00$.

Let $a=0, b=m\xi, c=n\xi, m,n\in\Z_+$ in (\ref{eq4.04}), we have
$$g(m\xi,n\xi)=\frac{\al m+\al n+mn\xi}\al\tau+\mu.$$
Substitute into (\ref{eq4.03}) with $a=m\xi, b=n\xi, c=k\xi, m,n,k\in\Z_+$,
then we obtain
$$\begin{aligned}
   0=&\left((\al n+\al k+nk\xi)\tau+\mu\al\right)
      \Huge(\left(\al m+\al(n+k+\zeta)+m(n+k+\zeta)\xi\right)\tau+\mu\al\Huge)\\
     &-\left((\al m+\al k+mk\xi)\tau+\mu\al\right)
      \Huge(\left(\al n+\al(m+k+\zeta)+n(m+k+\zeta)\xi\right)\tau+\mu\al\Huge),
\end{aligned}$$
considered as a polynomial in $m,n,k$,
whose leading term in the right hand side (if take $m\neq n$) is
$\tau^2\xi^2kmn(m-n)$.
This forces $\tau=0$.
So $\g\xi0=\mu$ for any $\xi\ing$.
Then taking $a=0$ in (\ref{eq4.04}) we show that $\g bc=\mu=\g00$ for all $b,c\ing$.

{\bf Case II}: $\al+\zeta\neq0$.
Then (\ref{eq4.17}) turns to
\begin{equation}\label{eq4.23}
 \g a0=\frac{(1+\theta\zeta)(\al+\zeta+a\al\theta)}{(\al+\zeta)(1+\theta(a+\zeta))}\g00.
\end{equation}
Set $a=0$ in (\ref{eq4.04}), then we get
\begin{equation}\label{eq4.24}
 \g bc=\frac1\zeta\left(f(b,c+\zeta)\g c0-\f bc\g{b+c}0\right).
\end{equation}
This shows that $g$ is a symmetric rational function on $G\times G$.
By Proposition \ref{prop4.4},
$g(a,b)=\frac{T(a)T(b)}{T(a+b+\zeta)}$
for some rational function $T(x)$ on $G$.
We note that $T(x)\neq0$ for all $x\ing$.

Put $b=c=0$ in (\ref{eq4.04}), we have
\begin{equation}\label{eq4.25}
 \frac{T(a)}{T(a+\zeta)}=\frac{T(0)}{T(\zeta)}\frac{\f a\zeta}{\al+\zeta}
  =\frac{T(0)}{T(\zeta)}\frac{(1+\theta\zeta)(\al+\zeta+a\al\theta)}{(\al+\zeta)(1+\theta(a+\zeta))}.
\end{equation}
Suppose $\theta\neq0$.
Then $\al\neq0$ by Proposition \ref{prop4.5}(1).
Hence $\frac{T(a)}{T(a+\zeta)}$ is not a constant.
Rewrite (\ref{eq4.25}),
$$\frac{T(a)}{T(a+\zeta)}=\frac{\al T(0)(1+\theta\zeta)}{T(\zeta)(\al+\zeta)}
  \frac{(a+\frac1\theta+\frac1{\al\theta}\zeta)}{(a+\frac1\theta+\zeta)}.$$
By Proposition \ref{prop4.5}(3) we see that
\begin{equation}\label{eq4.26}
 p=1-\frac1{\al\theta}\in\Z\setminus\{0\},\ \
 \frac{\al T(0)(1+\theta\zeta)}{T(\zeta)(\al+\zeta)}=1\text{ and }
 \frac{T(a)}{T(a+\zeta)}=\frac{(a+\frac1\theta+\frac1{\al\theta}\zeta)}{(a+\frac1\theta+\zeta)}.
\end{equation}

Notice that for any $a\ing$,
$g(a,-\zeta)=T(-\zeta)$ is a constant.
Set $c=-\zeta$ in (\ref{eq4.04}), we see
$$f(b,-\zeta)g(a,b-\zeta)=f(a,-\zeta)g(b,a-\zeta),$$
that is,
$$f(b,-\zeta)\frac{T(b-\zeta)}{T(b)}=f(a,-\zeta)\frac{T(a-\zeta)}{T(a)}\text{ for all }a,b\ing.$$
Take $b=\zeta$, and by (\ref{eq4.26}) we have
$$\frac{(a+\frac1\theta-\frac1{\al\theta}\zeta)(a+\frac1\theta+(\frac1{\al\theta}-1)\zeta)}
        {(a+\frac1\theta-\zeta)(a+\frac1\theta)}=
  \frac{(\frac1\theta-(\frac1{\al\theta}-1)\zeta)(\frac1\theta+\frac1{\al\theta}\zeta)}
        {\frac1\theta(\frac1\theta+\zeta)}
  \text{ for all } a\ing.$$
Hence the right hand side of the above equation equals to 1,
which implies $\al\theta=1$ and $p=0$, a contradiction to (\ref{eq4.26}).
So we have proved $\theta=0$ and $\f ab=\al+b$.
Then by (\ref{eq4.23}) and (\ref{eq4.24}) we see that
$\g ab=\g00$ for any $a,b\ing$.
}

Theorem \ref{thm4.1} follows directly from Proposition \ref{prop4.3}
and Proposition \ref{prop4.6}.
Thus by Theorem \ref{thm3.4}, \ref{thm3.4} and \ref{thm4.1}, we have
\begin{thm}\label{thm4.7}
Let $W$ be a CLSAS on $\W$ satisfying (\ref{eq4.01}) with arbitrary $f,g$.
Then $W$ must be isomorphic to one of $V_{\al,\theta}, V^{\gm,\lmd}$ or $W(\al,\mu,\zeta)$,
where $\al,\theta,\gm,\mu\in\C,\ \lmd,\zeta\ing$ such that $\theta=0$ or
$\theta^{-1}\not\ing,\ \lmd\neq\gm, \zeta\neq0$.
\end{thm}

A Novikov algebra $(N,\cdot)$
is a left-symmetric algebra with the multiplication $\cdot$ satisfying one more condition
$$(x\cdot y)\cdot z=(x\cdot z)\cdot y\text{ for any }x,y,z\ing.$$
Clearly, the left-symmetric algebra $W(\al,\mu,\zeta)$ is a Novikov algebra.
Moreover, we have
\begin{thm}
Any Novikov algebraic structure on $\W$ defined by (\ref{eq4.01}) must be isomorphic to
one of $W(\al,\mu,\zeta)$ with $\al, \mu\in\C$ and $\zeta\ing\setminus\{0\}$.
\end{thm}
\pf{
Notice by Theorem \ref{thm3.7} that $V_{\al,\frac1\al}\cong V_{\al,0}=W(\al,0,\zeta)$
if $\al\not\ing$.
Then by Theorem \ref{thm4.7} we only need to prove that
the left-symmetric algebras $V_{\al,\theta}$, $V^{\gm,\lmd}$ are not Novikov algebra
for any $\al,\theta,\gm\in\C,\lmd\ing$ such that $\theta\neq0,\al\theta\neq1,\theta^{-1}\not\ing$
and $\gm\neq\lmd$.

Suppose that $V_{\al,\theta}$ is a Novikov algebra for some
$\theta\neq0,\al\theta\neq1,\theta^{-1}\not\ing$.
From $(L_0L_b)L_c=(L_0L_c)L_b$, we derive that
$$\theta(\al\theta-1)bc(b-c)=0\text{ for any }b,c\ing,$$
which is a contradiction if we take $b,c\neq 0$ and $b\neq c$.
So such a $V_{\al,\theta}$ is not a Novikov algebra.

Suppose that $V^{\gm,\lmd}$ is a Novikov algebra.
Choose $a,b,c\ing\setminus\{0\}$ such that $\lmd+c\neq0, \lmd+a+c\neq0$ and $\lmd+a+b=0$.
From $(L_aL_b)L_c=(L_aL_c)L_b$ we get
$$ab(\lmd+c)=0,$$
which is a contradiction.
So $V^{\gm,\lmd}$ is not a Novikov algebra.
We complete the proof.}

The Balinskii-Novikov's construction \cite{BN} gives rise to a class of
infinite dimensional Lie algebras through some affinization of Novikov algebras.
Explicitly,
let $A$ be an algebra and define a binary operation $[\ ,\ ]$ on $A\otimes\C[t^{\pm1}]$ by
$$[x\otimes t^{i+1},y\otimes t^{j+1}]=
  (i+1)(xy)\otimes t^{i+j+1}-(j+1)(yx)\otimes t^{i+j+1},\ \ \forall x,y\in A, i,j\in\Z.$$
The Balinskii-Novikov's construction states that
$(A\otimes\C[t^{\pm1}],[\ ,\ ])$ is a Lie algebra if and only if $A$ is a Novikov algebra.

\begin{prop}
For any $\al,\mu\in\C, \zeta\ing\setminus\{0\}$,
there is an infinite dimensional Lie algebra induced from
the Novikov algebra $W(\al,\mu,\zeta)$,
which is isomorphic to the Lie algebra $\mathcal B(\al,\mu,\zeta)$
with a basis $\{x_{a,i}\mid a\ing,i\in\Z\}$ and Lie bracket
$$[x_{a,i},x_{b,j}]=\left((i+1)(\al+b)-(j+1)(\al+a)\right)x_{a+b,i+j}+
  \mu(i-j)x_{a+b+\zeta,i+j}.$$
\end{prop}
\pf{ The proposition is clear by setting $x_{a,i}=L_a\otimes t^{i+1}$.}

In particular, we get a class of Block-type Lie algebras (see \cite{B})
$\mathcal B(\al,0,\zeta)$, which are graded by a free group $G\oplus\Z$ of rank $\nu+1$.

\section{CLSAS on high rank Virasoro algebra $\V$}
\def\theequation{5.\arabic{equation}}
\setcounter{equation}{0}

In this section we study central extensions of the left-symmetric algebras
obtained in Section \ref{sec3},
whose sub-adjacent Lie algebra is $\V$.

Let $A$ be a left-symmetric algebra and $\phi:A\times A\longrightarrow\C$ a bilinear form satisfying
$$\phi(xy,z)-\phi(x,yz)=\phi(yx,z)-\phi(y,xz) \text{ for any }x,y,x\in A.$$
Denote $\hat A=A\oplus\C K$ and define
$$(x+aK)\ast(y+bK)=xy+\phi(x,y)K\text{ for }x,y\in A,\ a,b\in\C.$$
Then one can check that $(\hat A,\ast)$ is a left-symmetric algebra,
called the central extension of $A$ given by $\phi$.
Moreover, the bilinear form
$$\Phi(x,y)=\phi(x,y)-\phi(y,x)$$
defines a central extension of the sub-adjacent Lie algebra $L(A)$ of $A$.

Now apply these notions to a graded CLSAS on $\W$
and for convenience we write $\phif ab=\phi(L_a,L_b)$.
Then we have
\begin{align}
 &\phif ab-\phif ba=\frac 1{12}(b^3-b)\dt_{a+b,0};\label{eq5.01}\\
 &(b-a)\phif{a+b}c=\phif a{b+c}\f bc-\phif b{a+c}\f ac.\label{eq5.02}
\end{align}

Our main result in this section is the following
\begin{thm}\label{thm5.1}
A graded CLSAS $W$ on $\W$ has a central extension
whose sub-adjacent Lie algebra is $\V$ if and only if
$$W=V_{0,\theta}\text{ with }\theta\neq0 \text{ and }\theta^{-1}\not\ing.$$
Moreover, the central extension of $V_{0,\theta}$ is given by the bilinear form
$$\phi(L_a,L_b)=\frac1{24}\left(b^3-b-(\theta-\theta^{-1})b^2\right)\dt_{a+b,0}.$$
\end{thm}
\pf{
Let $W=V_{\al,\theta}$ or $V^{\gm,\lmd}$,
where $\al,\theta,\gm\in\C, \lmd\ing$ such that $\gm\neq\lmd$ and
$\theta=0$ or $\theta^{-1}\not\ing$.

Let $a=-b\neq0$ and $c=0$ in (\ref{eq5.02}) and by Lemma \ref{lem3.1} we have
$$\phif00=\frac1{24}\f00(b^2-1).$$
Since $\phif00$ is independent of $b$, it forces $\f00=0=\phif00$.
Then $\f b0=0,\f0b=b$ for any $b\ing$.
Set $b=c=0$ in (\ref{eq5.01}) and (\ref{eq5.02}). Then we get
$$\phif a0=\phif0a=0\text{ for any }a\ing.$$
Set $a=0$ in (\ref{eq5.02}) and we get
$$(b+c)\phif bc=0.$$
So we may write $\phif bc=\vf b\dt_{b+c,0}$ for some map $\varphi:G\longrightarrow\C$.
Let $b=-a$ in (\ref{eq5.01}) and we have
\begin{equation}\label{eq5.03}
 \vf a-\vf{-a}=\frac1{12}(-a^3+a).
\end{equation}
In the next we continue the proof in two cases.

{\bf Case I}: $W=V^{\gm,\lmd}$.
Notice that in this case $\lmd=\f00=0$, and we recall
$$\f bc=\begin{cases}
        c&\text{ if } b+c\neq0;\\
        \frac{c(\gm-c)}{\gm} &\text{ if }b+c=0.
  \end{cases}$$
Let $a+b+c=0$ in (\ref{eq5.02}) and we get
$$(a-b)\vf{a+b}=(a+b)(\vf a-\vf b).$$
Fix $0\neq c\ing$.
In the above equation, by taking $a=2c,b=-c$ we get
$$\vf{2c}=3\vf c+\vf{-c}=4\vf c+\frac1{12}(c^3-c),$$
and by taking $a=-2c,b=c$ we see
$$\vf{-2c}=3\vf{-c}+\vf c=4\vf c+\frac1{4}(c^3-c).$$
Put them into (\ref{eq5.01}) with $a=2c$ and we obtain
$$\frac16(-4c^3+c)=\vf{2c}-\vf{-2c}=-\frac16(c^3-c),$$
which implies $c=0$, a contradiction.
So the CLSAS $V^{\gm,\lmd}$ on $\W$ has no central extension required.

{\bf Case II}: $W=V_{\al,\theta}$.
Notice that in this case $\al=\f00=0$ and
$$\f bc=\frac{c(1+\theta c)}{1+\theta(b+c)}.$$
Let $a+b+c=0$ in (\ref{eq5.02}), then we have
$$\begin{aligned}
    &(b-a)\vf{a+b}=\vf a\f b{-a-b}-\vf b\f a{-a-b}\\
   =&-\frac{(a+b)\left(1-\theta(a+b)\right)}{1-\theta a}\vf a
     +\frac{(a+b)\left(1-\theta(a+b)\right)}{1-\theta b}\vf b.
  \end{aligned}$$
Denote $\etaf a=\frac{\vf a}{a(1-\theta a)}$.
Then from the above equation and (\ref{eq5.03}) we get
\begin{align}
 &(b-a)\etaf{a+b}=b\etaf b-a\etaf a;\label{eq5.04}\\
 &(1-\theta a)\etaf a+(1+\theta a)\eta(-a)=\frac1{12}(1-a^2).\label{eq5.05}
\end{align}
Now fix any $c\ing\setminus\{0\}$. Let $a=2c, b=-c$ in (\ref{eq5.04}), then we have
$3\etaf c=\eta(-c)+2\eta(2c)$.
Applying (\ref{eq5.05}) we get
\begin{equation}\label{eq5.06}
 \eta(2c)=\frac{\theta c+2}{\theta c+1}\etaf c+\frac{c^2-1}{24(\theta c+1)}.
\end{equation}
Notice here that $\theta c+1\neq0$ since $\theta=0$ or $\theta^{-1}\not\ing$.
Let $a=-2c, b=c$ in (\ref{eq5.04}) and we have
$3\eta(-c)=\eta(c)+2\eta(-2c).$
Also applying (\ref{eq5.05}) we get
\begin{equation}\label{eq5.07}
 \eta(-2c)=\frac{\theta c-2}{\theta c+1}\etaf c-\frac{c^2-1}{8(\theta c+1)}.
\end{equation}
Put (\ref{eq5.06}) and (\ref{eq5.07}) into (\ref{eq5.05}) with $a=2c$ and it gives
$$24\theta\etaf c-\theta-c=0.$$
This forces $\theta\neq0$ and $\etaf c=\frac1{24}(1+\theta^{-1}c)$.
So
$$\vf b=b\etaf b(1-\theta b)=\frac1{24}b(1-\theta b)(1+\theta^{-1}b),$$
and therefore
$$\phif ab=\frac1{24}\left(b^3-b-(\theta-\theta^{-1})b^2\right)\dt_{a+b,0}.$$
}

\begin{cor}
Any CLSAS on $\V$ with the multiplication
$$L_aL_b=\f abL_{a+b}+\phif abK,\ KL_a=L_aK=KK=0,$$
where $f,\phi$ are some functions on $G\times G$,
is isomorphic to the one defined by
$$L_aL_b=\frac{b(1+\theta b)}{1+\theta(a+b)}L_{a+b}
  +\frac K{24}\left(b^3-b-(\theta-\theta^{-1})b^2\right)\dt_{a+b,0},$$
where $\theta^{-1}\not\ing, \mathrm{Re}\theta>0$ or $\mathrm{Re}\theta=0, \mathrm{Im}\theta>0$.
\end{cor}

\end{document}